\documentclass{amsart}

\usepackage{amsmath,amssymb,amscd,amsfonts}

\newtheorem{theorem}{Theorem}
\newcommand{\bt}{\begin{theorem}}
\newcommand{\et}{\end{theorem}}
\newtheorem{lemma}{Lemma}
\newcommand{\bl}{\begin{lemma}}
\newcommand{\el}{\end{lemma}}
\newtheorem{corollary}{Corollary}
\newcommand{\bc}{\begin{corollary}}
\newcommand{\ec}{\end{corollary}}
\newcommand{\beq}{\begin{equation}}
\newcommand{\eeq}{\end{equation}}
\newcommand{\benum}{\begin{enumerate}}
\newcommand{\eenum}{\end{enumerate}}

\newcommand{\N}{\ensuremath{ \mathbf N }}

\DeclareMathOperator{\card}{\text{card}}

\newcommand{\bmat}{\left(\begin{matrix}}
\newcommand{\emat}{\end{matrix}\right)}
\DeclareMathOperator{\qand}{\quad\text{and}\quad}
\DeclareMathOperator{\qqand}{\qquad\text{and}\qquad}

\begin{document}
\title{The third positive element in the greedy $B_h$-set}
\author{Melvyn B. Nathanson}
\address{Department of Mathematics\\Lehman College (CUNY)\\Bronx, NY 10468} 
\email{melvyn.nathanson@lehman.cuny.edu}

\begin{abstract}
For $h \geq 1$, a $B_h$-set is a set of integers such that every integer $n$ has at most one representation 
in the form $n = a_{i_1} + \cdots + a_{i_h}$, where $a_{i_j} \in A$ for all $j = 1,\ldots, h$ 
and $a_{i_1} \leq \ldots \leq a_{i_h}$. 
The greedy $B_h$-set is the infinite set of nonnegative integers 
$\{a_0(h), a_1(h), a_2(h), \ldots \}$ 
constructed as follows:  If $a_0(h) = 0$ and $\{a_0(h), a_1(h), a_2(h),   \ldots, a_k(h) \}$ 
is a $B_h$-set, then $a_{k+1}(h)$ is the least positive integer such that 
$\{a_0(h), a_1(h), a_2(h), \ldots, a_k(h), a_{k+1}(h) \}$ is a $B_h$ set. 
One has $a_1(h) = 1$ and $a_2(h) = h+1$ for all $h$.  Elementary proofs are given that 
$a_3(h) = h^2+h+1$ for all $h \geq 1$ and that $a_k(h) \leq \sum_{i=0}^{k-1} h^i$ 
for all $h \geq 1$ and $k \geq 1$.
\end{abstract}



\maketitle

\section{Greedy Sidon sets and $B_h$-sets}

Let $A$ be a finite or infinite subset of the nonnegative integers $\N_0$ 
or, more generally, of any additive abelian semigroup $X$.   
For every positive integer $h$, the \emph{$h$-fold sumset} of $A$ is the set $hA$ consisting of all 
elements of $X$ that can be written as sums of $h$ not necessarily distinct elements of $A$. 
For example, if 
\[
A = \{0,1,3,7\} 
\]
then 
\[
2A = \{0,1,2,3,4,6,7,8,10,14\}
\]
and 
\[
3A = \{0,1,2,3,4,5,6,7,8,9,10,11,13,14,15,17,21\}.
\]

Let  
\[
n = a_{i_1} + a_{i_2} + \cdots + a_{i_h}
\]
and
\[
n = a_{j_1} + a_{j_2} + \cdots + a_{j_h}
\]
be two representations of $n$ as sums of $h$ elements of $A$.  These representations 
are \emph{equivalent} if there is a permutation $\sigma$ of $\{1,2, \ldots, h\}$ 
such that $a_{j_k} = a_{i_{\sigma(k)}}$ for all $k \in \{1,2, \ldots, h\}$. 
Inequivalent representations of $n$ are called \emph{distinct}. 
The \emph{representation function} $r_{h,A}(n)$ counts the number of 
distinct representations of $n$ as the sum of $h$ elements of $A$.  
The set $A$ is called a $B_h$-set if every element of $X$ has at most one distinct 
representation as the sum of $h$ elements of $A$, that is, if $r_{A,h}(x) = 0$ or $1$ 
for all $x \in X$ or, equivalently, if $r_{A,h}(x) = 1$ 
for all $x \in hA$.  
A $B_2$-set is also called a \emph{Sidon set}. 
 A finite set $A$ with $k$ elements is a $B_h$-set if and only if $\card(hA) = \binom{k+h-1}{h}$. 

We shall consider only sets of nonnegative integers.  
For a subset $A$ of $\N_0$, the representation function $r_{A,h}(n)$ 
counts the number of $h$-tuples 
$(a_{i_1}, a_{i_2}, \ldots, a_{i_h}) \in A^h$ 
such that 
\[
n = a_{i_1} + a_{i_2} + \cdots + a_{i_h}
\]
and
\[
a_{i_1} \leq  a_{i_2} \leq  \ldots \leq  a_{i_h}. 
\]

For example, if $A = \{0,1,3,7\}$, then $r_{2,A}(n) = 1$ for all $n \in 2A$ and so $A$ is a Sidon set, 
but $r_{3,A}(9) = 2$ (because $3 + 3 + 3 = 1 + 1 + 7$) and so $A$ is not a $B_3$-set. 
The infinite set $\{0, 1, 3, 3^2, 3^3,\ldots\}$ is a $B_2$-set 
but not a $B_3$-set (because $3^i + 3^i + 3^i = 3^{i+1} + 0 + 0$ for all $i = 0,1,2,\ldots$).

If $A$ is a $B_h$-set of integers, then the translation $t+A$ 
and the reflection $t-A$ are also $B_h$-sets.
It follows that if $A$ is a $B_h$-set of integers that is bounded below, 
then there is no loss of generality in assuming that $A$ is a set 
of nonnegative integers with $\min(A) = 0$.

The \emph{greedy $B_h$-set} is the infinite set $\{a_0(h), a_1(h), a_2(h), a_3(h), a_4(h), \ldots \}$ 
of nonnegative integers constructed inductively as follows: 
$a_0(h) = 0$ and, for $k \geq 1$, if  $\{a_0(h), a_1(h), \ldots, a_{k-1}(h)  \}$ is a $B_h$-set 
of integers with $a_0(h) < a_1(h) < \cdots <  a_{k-1(h)}$, 
then we choose $a_k(h)$ as the smallest positive integer such that $a_{k-1}(h) < a_k(h)$ and 
$\{a_0(h), a_1(h), \ldots, a_{k-1}(h), a_k(h) \}$ is a $B_h$-set. 
The greedy $B_1$-set is simply $\N_0 = \{0,1,2,3,4, \ldots\}$, that is, 
\[
a_k(1) = k \qquad \text{for all $k \in \N_0$.}  
\]
For all $h \geq 1$ we have 
\[
a_1(h) = 1 \qqand a_2(h) = h+1.  
\]

O'Bryant~\cite{obry23a} recently computed the 
initial  elements of greedy $B_h$-sets for small values of $h$.  
He observed that in every calculation the third positive element was $a_3(h) = h^2 + h + 1$.  
The purpose of this note is to prove that $a_3(h) = h^2 + h + 1$ for all $h \geq 1$.

Sidon sets and $B_h$-sets occur in  unexpected settings.  
For surveys and recent work on Sidon sets and $B_h$-sets, 
see 
Balogh, F\" uredi, and Roy~\cite{balo-fure-roy23}, 
Cilleruelo~\cite{cill17}, 
Halberstam and Roth~\cite{halb-roth66}, 
Nathanson~\cite{nath2022-202,nath2022-203}, 
O'Bryant~\cite{obry04,obry23}, 
and Pilatte~\cite{pila23}.

\section{The polynomial $a_3(h)$} 

\bl         \label{Bh-third:lemma:k+1}
Let $A$ be a finite set of nonnegative integers.  
If $A$ is a $B_h$-set with $\max(A) = a$, then $A \cup \{ha+1\}$ is a $B_h$-set.
If $a_k(h)$ is the $k$th positive element of the greedy $B_h$-set, then 
\[
a_{k+1}(h) \leq ha_k(h)+1. 
\]
\el

\begin{proof} 
Let $A = \{a_0 ,a_1, \ldots, a_k\}$ be a $B_h$-set with $ a_0 < a_1 < \cdots < a_k$, 
 let $a_{k+1} = ha_k+1$, and let $A' = A \cup \{a_{k+1} \}$.
If $A'$ is not a $B_h$-set, then there are  
distinct increasing sequences of nonnegative integers 
$0 \leq i_1 \leq i_2 \leq \cdots \leq i_h \leq k+1$ and 
$0 \leq j_1\leq j_2 \leq \cdots \leq  j_h \leq k+1$ 
such that 
\[
\sum_{t=1}^h a_{i_t} = \sum_{t=1}^h a_{j_t}.
\]
Cancelling summands that appear on both sides of this equation and renumbering, we obtain 
distinct  increasing sequences  of nonnegative integers 
 \[
0 \leq i_1 \leq i_2 \leq \cdots \leq i_{\ell} \leq k+1 \qqand 
0 \leq  j_1\leq j_2 \leq \cdots \leq  j_{\ell} \leq k + 1
 \]
 such that 
 \[
1 \leq \ell \leq h 
 \]
 and 
 \[
 \{ i_1,i_2, \ldots, i_{\ell} \} \cap \{  j_1, j_2, \ldots , j_{\ell} \} = \emptyset 
 \]
  and 
 \[
\sum_{t=1}^{\ell} a_{i_t} = \sum_{t=1}^{\ell} a_{j_t}.
\]
We must have $i_{\ell} = k+1$ or $j_{\ell} = k+1$ because $A$ is a $B_h$-set. 
If $i_{\ell} = k+1$, then $j_t \leq k$ for all $t \in \{1,\ldots, \ell\}$ and 
\[
ha_k+1 = a_{k+1} = a_{i_{\ell}} \leq \sum_{t=1}^{\ell} a_{i_t} = \sum_{t=1}^{\ell} a_{j_t} 
\leq \ell a_k \leq ha_k
\]
which is absurd.  The greedy algorithm immediately implies $a_{k+1}(h) \leq ha_k(h)+1$. 
This completes the proof. 
\end{proof}

\bc
The set $A = \{0\} \cup \left\{ \sum_{j=0}^{i-1} h^j : i  = 1,2,3,\ldots \right\}$ 
is an infinite $B_h$-set. 
\ec

\begin{proof}
Let $A_k = \{0\} \cup \left\{ \sum_{j=0}^{i-1} h^j : i  = 1,2,3,\ldots, k \right\}$ be a $B_h$-set.  
We have $\max(A_k) =   \sum_{j=0}^{k-1} h^j$ and $h\max(A_k)+1 =  \sum_{j=0}^{k} h^j.$ 
It follows from Lemma~\ref{Bh-third:lemma:k+1} that 
$A_{k+1} = \{0\} \cup \left\{ \sum_{j=0}^{i-1} h^j : i  = 1,2,3,\ldots, k+1 \right\}$ is a $B_h$-set. 
This completes the proof. 
\end{proof}

\bc                      \label{Bh-third:lemma:ak}
For all $h \geq 2$ and $k \geq 1$, the $k$th positive element $a_k(h)$ of the greedy $B_h$-set satisfies the inequality 
\[
a_k(h) \leq \sum_{i=0}^{k-1} h^{i} < h^{k-1} + 2h^{k-2}.
\]
\ec

\begin{proof}
Let $h \geq 2$.  
The proof is by induction on $k$.  For $k=1$ we have 
\[
a_1(h) = 1 = h^0 < h^0 +2h^{-1}.
\]  
Let $k \geq 1$ and let $\{0,a_1(h),\ldots, a_k(h)\}$ be the initial segment of the greedy $B_h$-set. 
If  $a_k(h) \leq  \sum_{i=0}^{k-1} h^i $, then, by Lemma~\ref{Bh-third:lemma:k+1},   
\begin{align*}
a_{k+1}(h) & \leq ha_k(h) +1 \leq  h\sum_{i=0}^{k-1} h^i + 1  = \sum_{i=0}^{k} h^i \\
& = h^k + h^{k-1} + \sum_{i=0}^{k-2} h^i 
 < h^k + 2h^{k-1}.
\end{align*}
This completes the proof. 
\end{proof}

\bt            \label{Bh-third:theorem}
Let  $h \geq 1$ and let 
\[
a_0(h) = 0 < a _1(h) = 1 < a_2(h) = h+1 < a_3(h)<\cdots
\]
 be the greedy $B_h$ set.  Then 
 \[
 a_3(h) = h^2 + h + 1.
 \]
\et

\begin{proof}
If $h = 1$, then $a_k(1) = k$ for all $k \geq 0$ and so $a_3(1) = 3 = 1^2 + 1 + 1$.

Let $h \geq 2$.
By Lemma~\ref{Bh-third:lemma:k+1}, 
\[
a_3(h) \leq ha_2(h)+1 =  h(h+1)+1 =  h^2+h+1.
\]
The proof  that $a_3(h) = h^2+h + 1$ is based on the observation that 
$a_3(h) \neq n$ if there exist nonnegative integers $x_0, x_1,x_2 $ and $y_0, y_1,y_2 $ 
such that 
\[
x_0 + x_1+x_2 = h-1
\]
\[
y_0 + y_1+y_2  = h
\]
and
\[
n = x_0a_0+x_1a_1+x_2a_2 = y_0a_0+y_1a_1+y_2a_2 
\]
or, equivalently, 
\[
n + x_1 + (h+1)x_2  = y_1 + (h+1)y_2. 
\]

Let $n = h^2 + h$.  We have 
\begin{align*}
n  & = n + \underbrace{0 + \cdots + 0 }_{\text{$h-1$ summands}} = 1 \cdot n + (h-1)a_0 \\
&  = \underbrace{ (h+1) + \cdots + (h+1)} _{\text{$h$ summands}} = ha_2
\end{align*}
and so $a_3(h) \neq n$.  

If $h+1 < n<  h^2+h$, then there exist unique integers $q$ and $r$ such that 
\[
n = q(h+1)+r  
\]
with 
\[
q \in \{1,2,\ldots, h-1\} \qand r \in \{0,1,\ldots, h\}.
\]

If $r = 0$, then 
\begin{align*}
n  
& = n + \underbrace{0 + \cdots + 0 }_{\text{$h-1$ summands}} = 1\cdot n + (h-1)a_0  \\
& =  q(h+1)  = \underbrace{ (h+1) + \cdots + (h+1)} _{\text{$q$ summands}}  
+ \underbrace{0+ \cdots +0} _{\text{$h-q$ summands}} \\
& =  (h-q)a_0 + qa_2  
\end{align*}
and so $a_3(h) \neq n$. 

If $r = 1$, then 
\begin{align*}
n  &  = n + \underbrace{0 + \cdots + 0 }_{\text{$h-1$ summands}} = 1\cdot n + (h-1)a_0  \\
& =   q(h+1) + 1 =  1 + \underbrace{ (h+1) + \cdots + (h+1)} _{\text{$q$ summands}}  + \underbrace{0+ \cdots +0} _{\text{$h-q -1$ summands}} \\
& = (h-q-1)a_0 + 1\cdot a_1 + qa_2  
\end{align*}
and so $a_3(h) \neq n$.

If $r \in \{2,3,\ldots, h\}$, then $1 \leq h+1-r \leq h-1$ and 
\begin{align*}
n + & (h+1-r) 
 = n +  \underbrace{1 + \cdots + 1} _{\text{$h+1-r$ summands}} + \underbrace{0 + \cdots + 0 }_{\text{$r-2$ summands}} \\ 
 & = 1\cdot n + (r-2)a_0 + (h+1-r)a_1  \\
& =  (q(h+1)+r) + (h+1-r) = (q+1)(h+1) \\
& =    \underbrace{ (h+1) + \cdots + (h+1)} _{\text{$q+1$ summands}}  
+ \underbrace{0+ \cdots +0} _{\text{$h-q -1$ summands}} \\
& =  (h-q-1)a_0  + (q+1)a_2 
\end{align*}
and so $a_3(h) \neq n$.  
Therefore, $a_3(h) \geq h^2 + h + 1$.
This completes the proof. 
\end{proof}

We have $a_4(h) \leq h^3 + h^2 + h + 1$ by Corollary~\ref{Bh-third:lemma:ak}.  
O'Bryant observed that his calculations suggest that $a_4(h)$ is the quasi-polynomial 
\[
a_4(h) = \begin{cases}
\frac{1}{2} (h+1)(h^2+2h+1) & \text{if $h$ is odd} \\
\frac{1}{2} (h+1)(h^2 + h + 2)& \text{if $h$ is even.}
\end{cases}
\]
Using the explicit value  $a_3(h) = h^2+h+1$ (Theorem~\ref{Bh-third:theorem}), 
Nathanson and O'Bryant~\cite{nath2023-221} have proved that this formula for $a_4(h)$ is correct.
Explicit formulae for $a_k(h)$ are not known for $k \geq 5$, 
nor is it not known if $a_k(h)$ is a quasi-polynomial of degree $k-1$ for $k \geq 5$. 

\def\cprime{$'$} \def\cprime{$'$} \def\cprime{$'$}
\providecommand{\bysame}{\leavevmode\hbox to3em{\hrulefill}\thinspace}
\providecommand{\MR}{\relax\ifhmode\unskip\space\fi MR }
\providecommand{\MRhref}[2]{%
  \href{http://www.ams.org/mathscinet-getitem?mr=#1}{#2}
}
\providecommand{\href}[2]{#2}


\begin{thebibliography}{1}

\bibitem{balo-fure-roy23}
J. Balogh, Z. F\"{u}redi, and S. Roy, \emph{An upper bound
  on the size of {S}idon sets}, Amer. Math. Monthly \textbf{130} (2023), no.~5,
  437--445.

\bibitem{cill17}
J. Cilleruelo, \emph{A greedy algorithm for {$B_h[g]$} sequences}, J.
  Combin. Theory Ser. A \textbf{150} (2017), 323--327.

\bibitem{pila23}
C. Pilatte, \emph{{A solution to the Erd\H os-S\' ark\" ozy-S\" os
  problem on asymptotic Sidon bases of order 3}}, {arXiv:2303:09659}, 2023.

\bibitem{halb-roth66}
H.~Halberstam and K.~F. Roth, \emph{{Sequences, Vol. 1}}, Oxford University
  Press, Oxford, 1966, Reprinted by Springer-Verlag, Heidelberg, in 1983.


\bibitem{nath2022-202}
M.~B. Nathanson,\emph{Sidon sets for linear forms}, J. Number Theory \textbf{239}
  (2022), 207--227. 





\bibitem{nath2022-203}
M.~B. Nathanson, \emph{The {B}ose-{C}howla argument for {S}idon sets}, J.
  Number Theory \textbf{238} (2022), 133--146. 
  
\bibitem{nath2023-221}
M.~B. Nathanson and K. O'Bryant, \emph{{The fourth element in a greedy
  $B_h$-set}}, preprint, 2023.

\bibitem{obry04}
K.~{O'Bryant}, \emph{A complete annotated bibliography of work related to
  {Sidon} sequences}, Electronic J. Combinatorics (2004), Dynamic Surveys DS
  11.

\bibitem{obry23}
K.~{O'Bryant}, \emph{Constructing thick {$B_h$} sets}, arXiv: 2308.12406, 2023.

\bibitem{obry23a}
K.~{O'Bryant}, Entries A365300--A365305 in The On-Line Encyclopedia of Integer Sequences, https://oeis.org/A365300--https://oeis.org/A365305.

\end{thebibliography}
\end{document}